\documentclass[a4paper,11pt]{amsart}

\usepackage{amssymb}
\usepackage{amsmath}
\usepackage{amsthm}

\newcommand{\1}{e}
\newcommand{\ii}{{\boldsymbol{i}}}

\newcommand{\CC}{\mathbb{C}}
\newcommand{\OO}{\mathcal{O}}
\newcommand{\PP}{\mathbb{P}}
\newcommand{\QQ}{\mathbb{Q}}
\newcommand{\RR}{\mathbb{R}}
\newcommand{\ZZ}{\mathbb{Z}}

\newcommand{\g}{\mathfrak{g}}
\newcommand{\h}{\mathfrak{h}}
\newcommand{\kk}{\mathfrak{k}}
\newcommand{\su}{\mathfrak{su}}
\newcommand{\ttt}{\mathfrak{t}}

\newcommand{\Aut}{\operatorname{Aut}}
\newcommand{\Ch}[1]{\Lambda_{#1}}
\newcommand{\diag}{\operatorname{diag}}
\newcommand{\Hom}{\operatorname{Hom}}
\newcommand{\Ho}[1]{\mathrm{H}^{#1}}
\newcommand{\Img}{\operatorname{Im}}
\newcommand{\Ker}{\operatorname{Ker}}
\newcommand{\ord}{\operatorname{ord}}
\newcommand{\rk}{\operatorname{rk}}
\newcommand{\Ru}[1]{#1_{\mathrm{u}}}
\newcommand{\tr}{\operatorname{tr}}

\newcommand{\Cone}{\mathcal{C}}
\newcommand{\ES}{E}
\newcommand{\Face}{\mathcal{U}}
\newcommand{\Fan}{\mathcal{F}}
\newcommand{\Kir}{\Psi}
\newcommand{\KKS}{\omega}
\newcommand{\Mom}{\Phi}
\newcommand{\Pol}{\mathcal{P}}
\newcommand{\sct}{s_\delta}
\newcommand{\Sph}{\mathrm{S}}
\newcommand{\sroot}{\sigma}
\newcommand{\Sym}{\mathrm{S}}
\newcommand{\Val}{\mathcal{V}}
\newcommand{\val}{\nu}
\newcommand{\vr}{\nu}
\newcommand{\WC}{\mathrm{C}^-}
\newcommand{\Wt}{\Lambda}
\newcommand{\wt}{\mu_\delta}
\newcommand{\Xbl}{\widetilde{X}}
\newcommand{\Xc}{\overline{X}}
\newcommand{\Xco}{\overline{X}^\circ}
\newcommand{\XG}{\mathfrak{X}}
\newcommand{\Xo}{X^\circ}
\newcommand{\Yc}{\overline{Y}}

\newcommand{\Mat}{\mathrm{Mat}}
\newcommand{\GL}{\mathrm{GL}}
\newcommand{\SL}{\mathrm{SL}}
\newcommand{\SO}{\mathrm{SO}}

\newcommand{\SU}{\mathrm{SU}}

\newcommand{\sst}{\text{\upshape ss}}

\renewcommand{\[}{\mathopen{[\mspace{-3mu}[}}
\renewcommand{\]}{\mathclose{]\mspace{-3mu}]}}
\renewcommand{\(}{\mathopen{(\mspace{-4mu}(}}
\renewcommand{\)}{\mathclose{)\mspace{-4mu})}}

\newtheorem{theorem}{Theorem}
\newtheorem*{SliceThm}{Slice Theorem}
\newtheorem*{LocStrThm}{Local Structure Theorem}
\newtheorem{proposition}[theorem]{Proposition}
\newtheorem{corollary}[theorem]{Corollary}

\theoremstyle{definition}
\newtheorem{example}[theorem]{Example}

\theoremstyle{remark}
\newtheorem{remark}[theorem]{Remark}

\title[On the orbit space of a maximal compact subgroup]
{On the orbit space \\ of a maximal compact subgroup \\ on a spherical homogeneous variety}

\author{Dmitry A. Timashev}

\address{Lomonosov Moscow State University, Faculty of Mechanics and Mathematics, Department of Higher Algebra, 119991 Moscow, Russia}
\address{Moscow Center of Fundamental and Applied Mathematics, 119991 Moscow, Russia}
\email{timashev@mccme.ru}



\begin{document}

\date{\today}

\begin{abstract}
Let $X=G/H$ be a spherical homogeneous variety for a complex reductive algebraic group~$G$. We prove that the orbit space of $X$ under the action of a maximal compact subgroup $K\subset G$ is homeomorphic to the valuation cone of~$X$. We also describe the relation between the orbit type stratification of the orbit space and the face stratification of the valuation cone.
\end{abstract}

\maketitle

\section{Introduction}

\subsection{}

Let $X=G/H$ be a homogeneous variety for a complex or real Lie group~$G$, $K\subset G$ be a maximal compact Lie subgroup, and $X/K$ be the orbit space for the action of $K$ on~$X$. It is well known from the theory of compact transformation groups \cite{compact} that $X/K$ is a Hausdorff topological space (in the quotient topology) stratified by smooth real manifolds according to the distribution of orbits by types corresponding to conjugacy classes of stabilizers. In this generality, nothing more substantial can be said about the orbit space, except that it is contractible. Contractibility stems from the following Karpelevich--Mostow theorem, which identifies the homogeneous variety $X$ with an equivariant vector bundle over a minimal $K$-orbit.

\begin{theorem}[\cite{fiber.Klein.I}, \cite{fiber.hom}, \cite{fiber.Klein.II}]\label{t:KM}
Suppose that the maximal compact subgroup $K$ is chosen (in its conjugacy class) in such a way that $K\cap H$ is a maximal compact subgroup in~$H$. Then there exists a real linear representation $K\cap H\to\GL(E)$ and a $K$-equivariant diffeomorphism $$X\simeq K\times^{K\cap H}E,$$ where the right-hand side of the formula is an equivariant vector bundle with the typical fiber $E$ over the compact homogeneous variety $K/(K\cap H)$.
\end{theorem}

\begin{corollary}
The orbit space $X/K$ is homeomorphic to $E/(K\cap H)$, hence contractible.
\end{corollary}

Thus $X/K$ is the topological cone over the orbit space $\Sph(E)/(K\cap H)$, where $\Sph(E)$ is the unit sphere in $E$ with respect to a $(K\cap H)$-invariant Euclidean metric. However the latter space $\Sph(E)/(K\cap H)$ can be arbitrarily complicated, because the representation of $K\cap H$ in $E$ is more or less arbitrary. Assuming $G$ be a complex linear algebraic group and $X$ be a complex homogeneous algebraic variety does not make the situation essentially simpler.

\subsection{}

In this paper we consider homogeneous varieties for linear algebraic groups over complex numbers. Specifically, we restrict ourselves to a quite special, though important, class of spherical homogeneous algebraic varieties. A (not necessarily homogeneous) algebraic variety $X$ is called \emph{spherical} if it is equipped with an action of a connected reductive algebraic group $G$ such that a Borel subgroup $B\subseteq G$ possesses a dense open orbit in~$X$. It follows that $G$ itself acts on $X$ with a dense open orbit, which is a spherical homogeneous variety.

The class of spherical varieties is remarkable in many aspects. It contains plenty of examples of algebraic varieties (both homogeneous and inhomogeneous) which are important in algebraic geometry, including Grassmannians, flag varieties, symmetric spaces, toric varieties, etc. For details see, e.g., \cite{sph} or \cite[Ch.~5]{hom-emb}.

The geometry of spherical varieties is controlled by several invariants of combinatorial and convex-geometric nature. One of these invariants is the \emph{valuation cone} $\Val_X$ of a spherical variety~$X$. It is a cosimplicial (i.e., defined by finitely many linearly independent linear inequalities) polyhedral convex cone in a real vector space, which is related to $G$-invariant valuations of the field of rational functions~$\CC(X)$; see a precise definition in~\ref{ss:sph}.

\subsection{}

Around 2016 Victor Batyrev conjectured that the orbit space $X/K$ for a spherical homogeneous variety $X$ is homeomorphic to~$\Val_X$. This conjecture, recorded only recently in~\cite{amoebae}, generalizes the classical Cartan decomposition in reductive groups and symmetric spaces; see Examples \ref{e:grp} and~\ref{e:sym}. Yet another argument in favour of Batyrev's conjecture is its ``non-Archimedean'' version describing the structure of the orbit space for the set $X\bigl(\CC\(\!\sqrt[\infty]t\,\)\bigr)$ of points of $X$ over the field of Puiseux series under the action of the group $G\bigl(\CC\[\!\sqrt[\infty]t\,\]\bigr)$ of points of $G$ over the ring of Puiseux series convergent at~$0$ (see~\ref{ss:non-Arch}).

Topologically, the cone $\Val_X$ either coincides with the ambient Euclidean space or is homeomorphic to a half-space in it. This dichotomy distinguishes between \emph{horospherical} (see Example~\ref{e:hor}) and other spherical homogeneous varieties. A particular case of Batyrev's conjecture is a criterion for horosphericity of a spherical homogeneous variety $X$ in terms of the topological structure of~$X/K$ \cite[Conj.\,7]{amoebae}.

However $\Val_X$ carries on a finer structure: it is a stratified manifold with corners, whose strata are the relative interiors of faces. A stronger version of Batyrev's conjecture \cite[Conj.\,6]{amoebae} states that $X/K$ is homomorphic to $\Val_X$ as a stratified space, so that the stratification of $X/K$ by orbit types corresponds to the stratification of $\Val_X$ by faces and the stabilizers for $K$-or\-bits are maximal compact subgroups in the so-called \emph{satellite subgroups} $H_{\Face}\subset G$ corresponding to faces $\Face\subseteq\Val_X$ (see~\ref{ss:sph}).

In this paper we prove the conjecture of Batyrev about homeomorphism between $X/K$ and $\Val_X$, and study a relation between the stratifications of these two spaces by orbit types and by faces, respectively. A feature of the proof is that a homeomorphism between $X/K$ and $\Val_X$ is constricted not directly, but via a third space, which is a polytope with several faces removed. This polytope $\Pol_{\Xc}$ is the image of the invariant moment map for the Hamiltonian action of $K$ on an equivariant projective compactification $\Xc$ of~$X$ (see~\ref{ss:mom}). It would be interesting to construct a natural homeomorphism between $X/K$ and $\Val_X$ directly.

We also prove that the stratification of $X/K$ by orbit types agrees with the stratification of $\Pol_{\Xc}$ by faces, meaning that the $K$-orbits corresponding to the interior points of a face of $\Pol_{\Xc}$ have one and the same type (Theorem~\ref{t:strat}). However the orbit types may coincide for different faces. Some general observations (Corollaries \ref{c:gen.sat} and~\ref{c:hor.sat}) and examples (see~\S\ref{s:ex}) justify the strong version of Batyrev's conjecture, but also show that the stratification of $X/K$ by orbit types may be coarser than the stratification of $\Val_X$ by faces (Example~\ref{e:SL_2^3/SL_2}). In general, the strong version of Batyrev's conjecture stays open.

\subsection{}

Now we describe the structure of the paper. In \S\ref{s:prelim} we introduce basic notation and provide necessary definitions and results on compact transformation groups and spherical varieties. In \S\ref{s:res} we state and prove our main results. Examples are considered in~\S\ref{s:ex}.

The author thanks the referee for careful reading of the initial version of the paper and for valuable comments, which helped to improve the presentation.

\section{Preliminaries}
\label{s:prelim}

\subsection{}

We keep the following basic notation, unless otherwise specified.

The multiplicative group of a field $F$ is denoted $F^\times=F\setminus\{0\}$.

We denote by $\1$ the unity element in any group~$P$. Given an action of $P$ on a set~$Y$, we denote by $Py$ the orbit of a point $y\in Y$ and by $P_y$ the stabilizer of $y$ in~$P$. The fixed point set for the action of $P$ on $Y$ is denoted~$Y^P$.

Let $P$ be a Lie group or an algebraic group, $Q\subseteq P$ be a closed Lie or algebraic subgroup, and $Z$ be a smooth manifold or a quasiprojective algebraic variety equipped with an action of~$Q$, respectively. Then there exists a unique, up to a $P$-equivariant isomorphism, smooth manifold, resp.\ algebraic variety, $Y$ acted on by $P$ which is equivariantly fibered over $P/Q$ in such a way that the fiber over the base point $\1Q$ is equivariantly isomorphic to~$Z$; see, e.g., \cite[Ch.\,2, \S3]{Lie-TG} and \cite[4.8]{inv.th}. We use the notation $$Y=P\times^QZ,\qquad\text{and}\qquad y=p*z$$ for the point $y\in Y$ in the fiber over~$pQ$ obtained from $z\in Z$ by the action of $p\in P$.

Let now $Y$ be any smooth manifold or algebraic variety acted on by $P$ which contains $Z$ as a $Q$-stable submanifold or subvariety. Then the embedding $Z\hookrightarrow Y$ extends in a unique way to a $P$-equivariant map $$P\times^QZ\to Y,\qquad p*z\mapsto y=p\cdot z,$$ given by applying elements of $P$ to points in~$Z$.

We denote the tangent space of a smooth manifold or algebraic variety $Y$ at a point $y\in Y$ as~$T_yY$.

The unipotent radical of a linear algebraic group $P$ is denoted by~$\Ru{P}$.

We use the following notation, unless otherwise specified:

\begin{itemize}
  \item $G$ is a connected reductive complex algebraic group;
  \item $B\subseteq G$ is a Borel subgroup;
  \item $T\subseteq B$ is a maximal algebraic torus;
  \item $K\subset G$ is a maximal compact subgroup;
  \item $X$ is a spherical homogeneous variety under an action of~$G$;
  \item $x_0\in X$ is a fixed base point;
  \item $H=G_{x_0}$ is the stabilizer in $G$ of the base point.
\end{itemize}
The choice of the base point defines an isomorphism $X\simeq G/H$. The subgroups $B$, $T$, $K$ are unique up to conjugation. We can (and will) choose them in a compatible way, so that $T\cap K$ will be a maximal (compact) tours in~$K$.

Lie groups and algebraic groups are denoted by capital Latin letters and their tangent Lie algebras are denoted by the same lowercase German letters. For instance, $\g$ is the Lie algebra of~$G$ and $\kk$ is the Lie algebra of~$K$.

\subsection{}\label{ss:comp}

We recall some results from the theory of compact transformation groups which we need. As a source, we use the book~\cite{compact}. In this subsection, we assume that $K$ is any compact Lie group and $X$ is a smooth real manifold equipped with a differentiable action of~$K$.

All $K$-orbits in $X$ are closed and compact. The orbit set $X/K$ is a Hausdorff space in the quotient topology induced by the natural map $X\to X/K$. The local structure of $X$ in a neighborhood of a given orbit is described by the following result.

\begin{SliceThm}
For any orbit $Kx\subset X$ there exists a submanifold $S\subset X$ of complementary dimension which is transversal to $Kx$ at~$x$, stable under the action of~$K_x$, and such that the map $$K\times^{K_x}S\to X,\qquad k*s\mapsto k\cdot s,$$ is a diffeomorphism onto a neighborhood $K\cdot S$ of~$Kx$.
\end{SliceThm}

The manifold $S$ is called a \emph{slice} for the action of $K$ on $X$ at the point~$x$. It can be constructed as the image of a sufficiently small neighborhood of $0$ in the normal space $N_xKx=(T_xKx)^\perp\subset T_xX$ with respect to a $K$-invariant Riemannian metric on~$X$ under the geodesic exponential map $T_xX\to X$.

The $K$-orbits in $X$ are distributed between classes of $K$-equivariant isomorphism, called \emph{orbit types}. An orbit type is determined by the conjugacy class of the stabilizer of a point in an orbit. Locally there are finitely many orbit types, and there is a natural hierarchy between them: the type of $Ky$ is higher than that of $Kx$ if, up to conjugation, $K_y\subset K_x$. It follows from the Slice Theorem that all orbits in a neighborhood of a given one have the same or higher type.

The division of orbits into types defines a stratification of~$X/K$. It follows from the Slice Theorem that the strata are locally closed in $X/K$ and carry on a structure of a smooth real manifold: the stratum containing the point corresponding to an orbit $Kx$ is diffeomorphic in a neighborhood of this point to the fixed point submanifold~$S^{K_x}$, which, in turn, is diffeomorphic to a neighborhood of $0$ in the vector space~$(N_xKx)^{K_x}$.

\subsection{}\label{ss:sph}

We introduce some basic notions and facts from the theory of spherical varieties. See, e.g., \cite{sph} or \cite{hom-emb} for more details.

We define a \emph{spherical variety} as an algebraic variety with an action of a connected reductive group~$G$ such that a Borel subgroup $B\subseteq G$ acts with an open orbit. The normality assumption on a variety is often included in the definition, but this issue is not essential for us, because we shall consider only smooth spherical varieties. The group $G$ itself also possesses an open orbit in~$X$, which is obviously a spherical homogeneous variety. We denote the open $B$-orbit in $X$ by~$\Xo$. A spherical $G$-variety contains finitely many $G$-orbits (and even $B$-orbits), and each $G$-orbit is a spherical homogeneous variety; see, e.g., \cite[Thm.\,2.1.2]{sph}, \cite[Cor.\,6.5]{hom-emb}.

The group $G$ acts on the set of rational functions on $X$ by the rule: $$(gf)(x)=f(g^{-1}x)\qquad (f\in\CC(X),\ g\in G,\ x\in X).$$ Consider \emph{$B$-semi-invariant} rational functions on~$X$, i.e., $f\in\CC(X)^\times$ such that $bf=\lambda(b)\cdot f$, $\forall b\in B$, for some character $\lambda:B\to\CC^{\times}$. These functions constitute a multiplicative subgroup $$\CC(X)^{(B)}\subseteq\CC(X)^\times.$$ A semi-invariant rational function $f=f_\lambda$ is determined by its weight $\lambda$ uniquely, up to a constant multiple. The weights of all $B$-semi-invariant rational functions constitute a sublattice $\Wt_X$ in the character lattice $\Ch{B}=\Ch{T}$ of~$B$ or, equivalently, of~$T$, which is called the \emph{weight lattice} of~$X$. There is an exact sequence
$$1 \longrightarrow \CC^\times \longrightarrow \CC(X)^{(B)} \longrightarrow \Wt_X \longrightarrow 0.$$

A \emph{$G$-invariant valuation} on $X$ is a $G$-invariant integer-valued discrete valuation of the function field $\CC(X)$ over the field of constants~$\CC$, i.e., a map ${\val:\CC(X)^\times\to\ZZ}$ satisfying the following conditions:
\begin{enumerate}
  \item $\val(f_1\cdot f_2)=\val(f_1)+\val(f_2)$, $\forall f_1,f_2\in\CC(X)^\times$;
  \item $\val(f_1+f_2)\ge\min\{\val(f_1),\val(f_2)\}$ whenever $f_1+f_2\ne0$;
  \item $\val(\CC^\times)=\{0\}$;
  \item $\val(gf)=\val(f)$, $\forall f\in\CC(X)^\times$, $g\in G$.
\end{enumerate}
If $X$ is embedded as an open $G$-orbit into a smooth spherical variety $\Xc$, then each boundary prime divisor $D\subseteq\Xc\setminus X$ defines a $G$-invariant valuation $\val=\val_D$, assigning to a rational function $f$ its vanishing order along~$D$. Conversely, any $G$-invariant valuation is of the shape $\val=k\cdot\val_D$ for some $\Xc$, $D$, and $k\in\ZZ_{\ge0}$ \cite[Prop.\,19.8, Prop.\,20.7]{hom-emb}.

Restricting a $G$-invariant valuation $\val$ to the subgroup $\CC(X)^{(B)}$, we can assign to $\val$ a vector $$\bar{\val}\in\Wt_X^*=\Hom(\Wt_X,\ZZ),\qquad\bar{\val}(\lambda)=\val(f_\lambda).$$
This vector determines $\val$ uniquely (D.~Luna and Th.~Vust, 1983; \cite[Cor.\,19.13]{hom-emb}). The vectors $\bar{\val}$ over all $G$-invariant valuations constitute the set $\Val_X(\ZZ)=\Val_X\cap\Wt_X^*$ of lattice points of a certain convex cone $\Val_X$ in the real vector space $\ES_X=\Hom(\Wt_X,\RR)$ spanned by the lattice~$\Wt_X^*$ (M.~Brion and F.~Pauer, 1987; F.~Knop, 1994; \cite[3.3.2]{sph}, \cite[Thm.\,20.3, Thm.\,21.1]{hom-emb}). The cone $\Val_X$ is called the \emph{valuation cone} of~$X$. It is a cosimplicial rational polyhedral cone of full dimension in~$\ES_X$, i.e., $$\Val_X=\{\vr\in\ES_X\mid\vr(\sroot_i)\le0,\ i=1,\dots,r\},$$ where $\sroot_1,\dots,\sroot_r$ is a collection (maybe empty) of linearly independent primitive vectors in the weight lattice~$\Wt_X$. The vectors $\sroot_i$, called \emph{spherical roots}, are non-negative linear combinations of positive roots of~$G$. They constitute a base of a certain root system in the vector space $\ES_X^*$ spanned by~$\Wt_X$, and $\Val_X$ is the negative (with respect to this base) Weyl chamber of this root system in~$\ES_X$. See, e.g., \cite[3.4.1]{sph}, \cite[22.4, 23.7]{hom-emb} for details.

\begin{remark}
Usually the valuation cone is defined in the \emph{rational} vector space $\ES_X(\QQ)=\Hom(\Wt_X,\QQ)$, but for us in this paper it will be more suitable to work with a real vector space. In our notation, the valuation cone in the traditional sense is just the set of rational points $\Val_X(\QQ)=\Val_X\cap\ES_X(\QQ)$.
\end{remark}

\begin{remark}
One can associate a vector $\bar{\val}\in\Wt_X^*$ to any (not necessarily $G$-in\-vari\-ant) integer-valued discrete valuation $\val$ of $\CC(X)$ over $\CC$, in particular, to the valuation $\val=\val_D$ defined by a prime divisor $D\subset X$. But in general, the vector $\bar{\val}$ no longer determines the valuation $\val$ uniquely.
\end{remark}

For any nonzero $G$-invariant valuation $\val$, consider a $G$-stable prime divisor $D$ on a spherical variety $\Xc$ in the complement of the open $G$-orbit $X$ such that $\val=k\cdot\val_D$, $k\in\ZZ_{>0}$. The normal bundle $N_{D,\Xc}$ is also a spherical $G$-variety \cite[Cor.\,2.4]{conormal}. The stabilizer $H_{\val}$ of a general point in the open $G$-orbit in $N_{D,\Xc}$ is called a \emph{satellite subgroup}. For $\val=0$, the satellite subgroup is, by definition, just $H\subset G$, the stabilizer of a point in~$X$. Up to conjugacy, the satellite subgroup $H_{\val}=H_\Face\subset G$ depends only on the face $\Face$ of~$\Val_X$ containing $\bar{\val}$ in its relative interior~\cite[Thm.\,1.1]{sat}. Satellite subgroups may be regarded as degenerations of~$H$.

\subsection{}\label{ss:non-Arch}

One can compute the $G$-invariant valuation~$\val=\val_D$ corresponding to a $G$-stable prime divisor $D\subset\Xc$ in the following way (see, e.g., \cite[\S24]{hom-emb}). Consider a smooth algebraic curve $C\subset\Xc$ intersecting $D$ transversally at a point~$y$ whose $G$-orbit is open in~$D$, so that $C\setminus\{y\}\subset X$. Any shifted curve $g\cdot C$, $g\in G$, has the same properties. If $f\in\CC(X)$ is any rational function, then for general $g\in G$ the point $gy$ does not belong to components of the divisor of $f$ different from~$D$. Thus the valuation can be computed by the formula
\begin{equation}\label{eq:val}
\val(f)=\val_{gy}\bigl(f|_{gC}\bigr)=\ord_{t=0}f\bigl(gx(t)\bigr),
\end{equation}
where $g\in G$ is a general element and $t$ is a uniformizing variable parameterizing points $x(t)\in C$ in a neighborhood of $y=x(0)$.

Replacing an actual neighborhood of $y$ by a formal one, one may assume that $x(t)$ is a point of $X$ over the field of Laurent series $\CC\(t\)$ in the formal variable~$t$. Conversely, any point $x(t)\in X\bigl(\CC\(t\)\bigr)$ defines a $G$-invariant valuation $\val=\val_{x(t)}$ by the formula~\eqref{eq:val}. The group $G\bigl(\CC\(t\)\bigr)$ acts on $X\bigl(\CC\(t\)\bigr)$ and the subgroup $G\bigl(\CC\[t\]\bigr)$ preserves~$\val$. The Galois group $\Gamma=\Aut\CC\[t\]$ consisting of changes of the uniformizing parameter also acts on $X\bigl(\CC\(t\)\bigr)$ and preserves~$\val$.

Replacing $\val$ with $k\cdot\val$ ($k\in\ZZ_{>0}$) is equivalent to replacing $t$ with~$t^k$. One may consider not only integer but also rational~$k$. Then one has to replace the field of Laurent series $\CC\(t\)$ with its algebraic closure, the field of Puiseux series $\CC\(\!\sqrt[\infty]t\,\)$, the ring of formal power series $\CC\[t\]$ with the ring $\CC\[\!\sqrt[\infty]t\,\]$ of Puiseux series convergent at $t=0$, and $\Gamma$ with the group $\Gamma_\infty$ of algebra automorphisms of $\CC\[\!\sqrt[\infty]t\,\]$ preserving the order at $0$. Observe that the group $G\bigl(\CC\(\!\sqrt[\infty]t\,\)\bigr)$ acts on the set $X\bigl(\CC\(\!\sqrt[\infty]t\,\)\bigr)$ transitively, with stabilizer $H\bigl(\CC\(\!\sqrt[\infty]t\,\)\bigr)$.

The following result of Luna and Vust can be viewed as a ``non-Archimedean analogue'' of Batyrev's conjecture.

\begin{theorem}[{\cite[\S4]{LV}, \cite[\S24]{hom-emb}}]
There are bijections
\begin{align*}
X\bigl(\CC\(t\)\bigr)/\bigl(G(\CC\[t\])\rtimes\Gamma\bigr)                                   &\simeq \Val_X(\ZZ),\\ X\bigl(\CC\(\!\sqrt[\infty]t\,\)\bigr)/\bigl(G(\CC\[\!\sqrt[\infty]t\,\])\rtimes\Gamma_\infty\bigr) &\simeq \Val_X(\QQ).
\end{align*}
\end{theorem}

The subgroups $G\bigl(\CC\[t\]\bigr)\subset G\bigl(\CC\(t\)\bigr)$ and $G\bigl(\CC\[\!\sqrt[\infty]t\,\]\bigr)\subset G\bigl(\CC\(\!\sqrt[\infty]t\,\)\bigr)$ are analogues of the maximal compact subgroup $K\subset G$. If $G$ is defined over $\ZZ$, then after replacing $\CC$ with a finite field $F$ we get that $G\bigl(F\[t\]\bigr)$ is indeed a maximal compact subgroup of $G\bigl(F\(t\)\bigr)$, the group of points of $G$ over the non-Archimedean local field $F\(t\)$, in the $t$-adic topology.

\subsection{}\label{ss:mom}

Let $\Xc$ be a smooth spherical variety with an open $G$-orbit~$X$. Since the open $B$-orbit $\Xo\subseteq X\subseteq\Xc$ is an affine variety, its boundary $\Xc\setminus\Xo$ consists of finitely many $B$-stable prime divisors $D_1,\dots,D_l$. Among them, there may occur both $G$-stable divisors lying in the complement of $X$ and $B$-stable divisors intersecting~$X$. Let $\val_i=\val_{D_i}$ be the respective valuations of $\CC(X)$ and $\bar{\val}_i\in\Wt_X^*$ be the respective lattice vectors in~$E_X$.

A $G$-orbit $Y\subset\Xc$ is uniquely determined by the set of $B$-stable divisors $D_{i_1},\dots,D_{i_k}$ containing~$Y$ \cite[Thm.\,3.1.4]{sph}, \cite[Prop.\,14.1]{hom-emb}. The set of respective vectors $\bar{\val}_{i_1},\dots,\bar{\val}_{i_k}$ can be included into a basis of the lattice~$\Wt_X^*$ \cite[Cor.\,4.3.19]{sph}. The subvariety $\Yc$ is the center of a $G$-invariant valuation $\val$ if and only if $\bar{\val}$ belongs to the relative interior of the cone $\Cone_Y\subset\ES_X$ spanned by $\bar{\val}_{i_1},\dots,\bar{\val}_{i_k}$; in particular, the relative interior of $\Cone_Y$ intersects~$\Val_X$ \cite[Thm.\,3.1.7]{sph}, \cite[15.1]{hom-emb}.

The variety $\Xc$ is called \emph{toroidal} if for any $G$-orbit $Y$ all the divisors $D_i$ containing $Y$ are $G$-stable. The respective cones $\Cone_Y$ constitute a smooth simplicial fan~$\Fan_{\Xc}$ contained in~$\Val_X$. (Recall that a \emph{smooth simplicial fan} in a vector space spanned by a lattice is a finite collection of convex cones such that each cone is generated by a subset of a basis of the lattice, any two different cones intersect in their common face, and all faces of each cone in the fan also belong to the fan.) Also, for any two orbits $Y,Y'\subset\Xc$ we have $\Yc\supset Y'$ if and only if $\Cone_Y$ is a face of~$\Cone_{Y'}$. The variety $\Xc$ is complete if and only if its fan $\Fan_{\Xc}$ covers the whole valuation cone~$\Val_X$. One can obtain a toroidal variety from any spherical variety by a sequence of blow-ups of orbit closures. See \cite[3.3]{sph}, \cite[\S29]{hom-emb} for details.

If $\Xc$ is projective, then there exists a $G$-equivariant embedding ${\Xc\hookrightarrow\PP(V)}$, where $V$ is a finite-dimensional rational $G$-module. The hyperplane section of $\Xc$ defined by a $B$-semi-invariant linear form (i.e., a highest weight vector in~$V^*$) is a $B$-stable very ample effective divisor of the shape $$\delta=m_1D_1+\dots+m_lD_l.$$ Without loss of generality we may assume that $V=\Ho0\bigl(\Xc,\OO(\delta)\bigr)^*$. Let $\sct\in\Ho0\bigl(\Xc,\OO(\delta)\bigr)$ be a $B$-semi-invariant global section of the line bundle $\OO(\delta)$ corresponding to~$\delta$ and $\wt\in\Ch{B}$ be its eigenweight.

Consider a convex polytope $\Pol_{\Xc}$ in the affine space $\ES_X^*+\wt$ (where $\ES_X^*=\Wt_X\otimes\RR$) defined by the formula
$$\Pol_{\Xc}=\bigl\{\mu\in\ES_X^*+\wt\bigm|\bar{\val}_i(\mu-\wt)\ge-m_i,\ i=1,\dots,l\bigr\}.$$
Replacing $\delta$ with $n\delta$ yields a dilation of $\Pol_{\Xc}$ with scale factor~$n$.

The linear representation of $G$ in the space $\Ho0\bigl(\Xc,\OO(\delta)\bigr)$ is \emph{multiplicity free}, i.e., any irreducible $G$-module occurs in its decomposition at most once; see, e.g., \cite[Thm.\,2.1.2]{sph}. The highest weight vectors of the irreducible submodules in $\Ho0\bigl(\Xc,\OO(\delta)\bigr)$, i.e., the $B$-semi-invariant global sections of $\OO(\delta)$ are of the shape $s=f_\lambda\sct$, where $\lambda\in(\Pol_{\Xc}-\wt)\cap\Wt_X$. Thus the highest weights $\mu=\lambda+\wt$ of the irreducible submodules in $\Ho0\bigl(\Xc,\OO(\delta)\bigr)$ are exactly the lattice points of~$\Pol_{\Xc}$.

A description of the space of $K$-orbits in $\Xc$ in terms of symplectic geometry was suggested by Brion~\cite{moment}. We reproduce it here.

There exists a $K$-invariant Hermitian inner product $(\,{\cdot}\,|\,{\cdot}\,)$ on~$V$. It defines a K\"ahler metric on~$\Xc$, known as the Fubini--Study metric. The imaginary part of the Fubini--Study metric is a $K$-invariant real symplectic form. Thus $\Xc$ comes equipped with a structure of a real symplectic manifold with a Hamiltonian action of~$K$. The respective moment map $\Mom:\Xc\to\kk^*$ is given by the formula
\begin{equation}\label{eq:mom}
\langle\Mom(x),\xi\rangle=\frac1\ii\cdot\frac{(v|\xi{v})}{(v|v)},\qquad\forall\xi\in\kk,
\end{equation}
where $v\in V$ is a vector lying over $x\in\PP(V)$. Identifying $\kk^*$ with $\kk$ by means of a $K$-invariant Euclidean inner product on~$\kk$, one may regard $\Mom$ as a map to $\kk$ and interpret $\langle\,{\cdot}\,,\,{\cdot}\,\rangle$ in \eqref{eq:mom} as this inner product.

The lattice $\Ch{B}=\Ch{T}$ embeds into the space $\ttt^*$ (by assigning to a character its differential at~$\1$) and spans its real form $\ii(\ttt\cap\kk)^*$ consisting of the linear functions that take pure imaginary values on~$\ttt\cap\kk$. Denote by $\WC$ the negative Weyl chamber in~$\ttt\cap\kk$, which consists of the vectors on which all positive roots take values in~$\ii\RR_{\le0}$. The chamber $\WC$ is a fundamental domain for the adjoint action of $K$ on~$\kk$: each $K$-orbit intersects $\WC$ in a single point. Under identification of $\ii(\ttt\cap\kk)^*$ with $\ttt\cap\kk$ by means of the invariant inner product, the set of lattice points $\WC\cap\Ch{T}$ is exactly the set of \emph{antidominant weights}, which are opposite to the highest weights of the irreducible representations of~$G$.

Composing the moment map $\Mom$ with the quotient map $\kk\to\WC$ for the adjoint action of $K$ yields the \emph{invariant moment map} $\Kir:\Xc\to\WC$, also called the \emph{Kirwan map} (in honor of F.~Kirwan).

\begin{theorem}[{\cite[\S5]{moment}}]\label{t:mom}
The fibers of the Kirwan map are the $K$-orbits on $\Xc$ and the image of this map is the polytope $-\Pol_{\Xc}\subset\WC$. In other words, $\Kir$~is the quotient map for the action of $K$ on~$\Xc$, and $$\Img\Kir=-\Pol_{\Xc}\simeq\Xc/K.$$
\end{theorem}

In view of this result, $\Pol_{\Xc}$ is sometimes called the \emph{moment polytope} of a projective spherical variety~$\Xc$.

\section{Main results}
\label{s:res}

\subsection{}

Here we prove the conjecture of Batyrev on the topological structure of the space of $K$-orbits.

\begin{theorem}\label{t:coarse}
For a spherical homogeneous $G$-variety $X$, the orbit space $X/K$ of a maximal compact subgroup $K\subset G$ is homeomorphic to the valuation cone~$\Val_X$.
\end{theorem}

Recall that a homogeneous variety $X=G/H$ is called \emph{horospherical} if $H$ contains a maximal unipotent subgroup of~$G$. Horospherical varieties are spherical (cf.~Example~\ref{e:hor}). It is well known (see, e.g., \cite[21.3]{hom-emb}) that a spherical homogeneous variety $X$ is horospherical if and only if $\Val_X=\ES_X$.

\begin{corollary}[{\cite[Conj.\,7]{amoebae}}]
$X$ is horospherical if and only if $X/K$ is homeomorphic to a Euclidean space.
\end{corollary}

In order to prove Theorem~\ref{t:coarse}, we choose a smooth projective toroidal compactification $\Xc$ of $X$ determined by a smooth simplicial fan $\Fan_{\Xc}$ subdividing the valuation cone~$\Val_X$; see, e.g., \cite[29.2]{hom-emb}. ($\Xc$~is complete, since $\Fan_{\Xc}$ covers the whole~$\Val_X$, and projectivity can be obtained by further subdivision in order to get a strictly convex piecewise linear function on~$\Fan_X$ \cite[Cor.\,3.2.12]{sph}.) We embed $\Xc$ equivariantly into a projective space~$\PP(V)$, as in~\ref{ss:mom}. The following two results are contained in \cite[5.3]{var.sph}; we prove them in our situation for the reader's convenience.

\begin{proposition}\label{p:norm.fan}
The fan $\Fan_{\Xc}$ is a part of the normal fan of the polytope~$\Pol_{\Xc}$.
\end{proposition}

\begin{proof}
Since the hyperplane section divisor $\delta=\sum_im_iD_i$ is ample on~$\Xc$, for any $G$-orbit $Y\subset\Xc$ there is a point $\mu=\mu_Y\in\Pol_{\Xc}$ such that the defining inequalities
\begin{equation}\label{eq:pol}
\bar{\val}_i(\mu-\wt)\ge-m_i,\qquad i=1,\dots,l,
\end{equation}
at this point turn into equalities exactly when $D_i\supseteq Y$; see \cite[Thm.\,3.2.9]{sph}, \cite[Thm.\,17.18]{hom-emb}. Indeed, for $n$ sufficiently big one can find a $B$-semi-in\-vari\-ant global section $s\in\Ho0\bigl(\Xc,\OO(n\delta)\bigr)$ such that $s\ne0$ on $Y$ and $s=0$ on any $D_i\not\supseteq Y$, and then put $\mu_Y=\tilde\mu/n$, where $\tilde\mu$ is the eigenweight of~$s$. If the orbit $Y$ is closed, i.e., the cone $\Cone_Y$ is full-dimensional, then the point $\mu_Y\in\Wt_X+\wt$ is unique and is a vertex of~$\Pol_{\Xc}$. In general, such points fill the relative interior of the face of~$\Pol_{\Xc}$ distinguished by the equations
\begin{equation}\label{eq:face}
\bar{\val}_i(\mu-\wt)=-m_i,\qquad\forall D_i\supseteq Y,
\end{equation}
with the normal cone~$\Cone_Y$.
\end{proof}

The closure $\Yc$ of any $G$-orbit $Y\subset\Xc$ is also a smooth projective spherical variety, because $\Yc$ is the intersection of the prime $G$-stable divisors $D_i$ containing~$Y$, which constitute a divisor with normal crossings; see, e.g., \cite[Thm.\,29.2]{hom-emb}. Let us describe the moment polytope of~$\Yc$.

\begin{proposition}\label{p:mom.orb}
$\Pol_{\Yc}$ is the face of $\Pol_{\Xc}$ defined by the equations~\eqref{eq:face}, with the normal cone~$\Cone_Y$.
\end{proposition}

\begin{proof}
The restriction map $$\Ho0\bigl(\Xc,\OO(n\delta)\bigr)\longrightarrow\Ho0\bigl(\Yc,\OO(n\delta)\bigr)$$ is surjective for big (in fact, any)~$n$. Also, by complete reducibility of linear representations of~$G$, $B$-semi-invariant sections on $\Yc$ can be extended to $B$-semi-invariant sections on~$\Xc$. A $B$-semi-invariant section $s\in\Ho0\bigl(\Xc,\OO(n\delta)\bigr)$ of weight $\tilde\mu$ has nonzero restriction to $\Yc$ if and only if $\mu=\tilde\mu/n\in\Pol_{\Xc}$ satisfies~\eqref{eq:face}. It follows that the polytope $\Pol_{\Yc}$ and the face of $\Pol_{\Xc}$ defined by~\eqref{eq:face} have one and the same set of rational point, hence coincide.
\end{proof}

\begin{proof}[Proof of Theorem~\ref{t:coarse}]
Theorem \ref{t:mom} yields a homeomorphism
\begin{equation}\label{eq:X/K}
X/K\simeq\Pol_{\Xc}\setminus\bigcup_Y\Pol_{\Yc},
\end{equation}
where the union is taken over all $G$-orbits $Y\subseteq\Xc\setminus X$. It follows from Propositions \ref{p:norm.fan} and \ref{p:mom.orb} that $\bigcup_Y\Pol_{\Yc}$ is the union of the faces of $\Pol_{\Xc}$ whose normal cones constitute the fan $\Fan_{\Xc}$ covering the valuation cone~$\Val_X$.

If $X$ is horospherical, then $\Val_X=\ES_X$ and $\bigcup_Y\Pol_{\Yc}$ is the whole boundary of~$\Pol_{\Xc}$, which is homeomorphic to the boundary sphere of a ball in~$\ES_X^*$. Otherwise $\Val_X$ has boundary and is homeomorphic to a half-space in~$\ES_X$ and $\bigcup_Y\Pol_{\Yc}$ is a part of the boundary of~$\Pol_{\Xc}$, which is homeomorphic to a closed hemisphere in the boundary sphere. In both cases, ${\Pol_{\Xc}\setminus\bigcup_Y\Pol_{\Yc}}$ is homeomorphic to~$\Val_X$.
\end{proof}

\subsection{}

In order to study the stratification of $X/K$, we have to examine properties of the moment map ${\Mom:\Xc\to\kk^*}$ and the Kirwan map $\Kir:\Xc\to\WC$ in more detail.

Let $V$ be a finite-dimensional rational $G$-module. Recall (see, e.g., \cite[4.6]{inv.th}) that a point $x\in\PP(V)$ corresponding to a nonzero vector $v\in V$ is called \emph{semistable} if there exists a homogeneous invariant polynomial $F\in\CC[V]^G$ of positive degree such that $F(v)\ne0$. The set of semistable points $\PP(V)^{\sst}$ is open in~$\PP(V)$. The following result refines (by adding information on stabilizers) a criterion of closedness of orbits, due to G.\,Kempf and L.~Ness, which is well known in Invariant Theory.

\begin{proposition}\label{p:KN-crit}
Let $V$ be a rational $G$-module equipped with a $K$-invariant Hermitian inner product, $x\in\PP(V)$, and $v\in V$ be a nonzero vector lying over~$x$. Then:
\begin{enumerate}
  \item\label{i:closed} the orbit $Gv$ is closed in~$V$ and $Gx$ is closed in $\PP(V)^{\sst}$ if and only if $\Mom(x')=0$ for some $x'\in Gx$;
  \item\label{i:crit} all such points $x'$ and all overlying vectors $v'\in Gv$, respectively, constitute a single $K$-orbit;
  \item\label{i:stab} their stabilizers $K_{x'}$ and $K_{v'}$ are maximal compact subgroups in the reductive groups $G_{x'}$ and~$G_{v'}$, respectively.
\end{enumerate}
\end{proposition}

\begin{proof}
A proof of \eqref{i:closed} and \eqref{i:crit} can be found in \cite[6.12]{inv.th}, see also \cite[Rem.\,7.5, Thm.\,8.10]{Kir}. We have to recall these arguments in part, in order to prove~\eqref{i:stab}.

First note that $\Mom(x')=0$ if and only if $v'$ is orthogonal to the tangent space $T_{v'}Gv=\g{v'}$. Such vectors $v'$ are the critical points of the norm-square function $(v'|v')$ on the orbit~$Gv$. Thus the existence of such vectors is necessary for the orbit to be closed.

We may assume further without loss of generality that ${\Mom(x)=0}$.

The polar decomposition of the reductive group $G=K\cdot\exp(\ii\kk)$ \cite[1.6.4]{Lie-3} allows us to express any $g\in G$ uniquely as $g=k\cdot\exp\xi$, where $k\in K$, $\xi\in\ii\kk$. Put $g(t)=k\cdot\exp(t\xi)$ ($t\in\RR$) and consider a decomposition $v=v_1+\dots+v_m$ into a sum of eigenvectors of the Hermitian operator $\xi$ with eigenvalues $\lambda_1,\dots,\lambda_m\in\RR$.

We have $(v|\xi{v})=\sum_i\lambda_i(v_i|v_i)=0$. Therefore, either $\xi{v}=0$, or there occur numbers of different signs among~$\lambda_i$'s. In the first case, $\exp\xi\in G_v$ and ${gv\in Kv}$. In the second case, the function $$q(t)=\bigl(g(t)v\!\bigm|\!g(t)v\bigr)=\sum_ie^{2\lambda_it}(v_i|v_i)$$ is strictly convex and tends to $+\infty$ as $t\to\pm\infty$, whence has a unique critical point $t=0$, which is the minimum point.

It follows that all critical points of the norm-square function on $Gv$ are minimum point and constitute a single orbit $Kv$ of the group~$K$. Besides, $gv=v$ if and only if $kv=v$, $\xi{v}=0$. Therefore $K_v$ is a maximal compact subgroup in~$G_v$.

Since $G_v$ is a subgroup of finite index in~$G_x$, we have $gx=x$ if and only if $gv=\omega\cdot{v}$ for some root of unity~$\omega$, which is equivalent to $kv=\omega\cdot{v}$, $\xi{v}=0$, i.e., $kx=x$. Hence $K_x$ is a maximal compact subgroup in~$G_x$.
\end{proof}

It follows from Theorem \ref{t:mom} that for each $\mu\in\Pol_{\Xc}$ the fiber of the moment map $\Mom^{-1}(-\mu)$ is an orbit of the group~$K_{\mu}$, the stabilizer of~$\mu$ (and~$-\mu$) in the (co)adjoint representation. Let us describe it more explicitly.

Consider all the $B$-stable divisors $D_i$ such that the inequalities \eqref{eq:pol} turn into equalities. These equalities define a face $\Pol(\mu)$ of~$\Pol_{\Xc}$ whose relative interior contains~$\mu$. The complement of the union of the remaining divisors $D_i$ is an affine open subset $\Xc(\mu)\subset\Xc$. If $\mu$ is a rational point in~$\Pol(\mu)$, then for sufficiently big $n$ there exists a $B$-stable section $s\in\Ho0\bigl(\Xc,\OO(n\delta)\bigr)$ of eigenweight $\tilde\mu=n\mu$, and $\Xc(\mu)=\Xc_s$ is the respective principal open subset.

The parabolic subgroup $P(\mu)\subseteq G$ containing $B$ and consisting of the elements preserving $\Xc(\mu)$ or, equivalently, all the divisors $D_i$ in the complement of~$\Xc(\mu)$ depends only on the face $\WC(-\mu)$ of the Weyl chamber~$\WC$ containing $-\mu$ in its relative interior. If $\mu$ is a rational point, then, in the above notation, $P(\mu)$~is the stabilizer of the line spanned by~$s$. The roots of $P(\mu)$ are those roots $\alpha$ of~$G$ which satisfy $\langle\alpha,\mu\rangle\ge0$. There is a Levi decomposition $P(\mu)=\Ru{P(\mu)}\rtimes L(\mu)$, where $\Ru{P(\mu)}$ is the unipotent radical and $L(\mu)$ is the standard Levi subgroup of $P(\mu)$ containing~$T$, whose roots are the roots of $G$ orthogonal to~$\mu$. Note that $K_\mu$ is a maximal compact subgroup in~$L(\mu)$.

The following result is a particular case of the so-called \emph{local structure theorem} of Brion--Luna--Vust; see, e.g., \cite[Thm.\,2.3.2]{sph} or \cite[Thm.\,4.6]{hom-emb}.

\begin{LocStrThm}
There exists an $L(\mu)$-stable closed subvariety $Z(\mu)\subseteq\Xc(\mu)$ such that
$$\Xc(\mu)\simeq P(\mu)\times^{L(\mu)}Z(\mu)\simeq\Ru{P(\mu)}\times Z(\mu).$$
\end{LocStrThm}

The variety $Z(\mu)$ is called a \emph{BLV-slice} of the affine chart $\Xc(\mu)$ of the spherical variety~$\Xc$. It is an affine spherical $L(\mu)$-variety and therefore contains a dense open $L(\mu)$-orbit, and also contains a unique closed $L(\mu)$-or\-bit.

\begin{proposition}\label{p:fib.mom}
$\Mom^{-1}(-\mu)=K_{\mu}\cdot{x}\subseteq\Xc(\mu)$ is a single $K_{\mu}$-orbit, and $P(\mu)x$ is the unique closed $P(\mu)$-orbit in~$\Xc(\mu)$. It intersects $Z(\mu)$ in the unique closed orbit $L(\mu)z$, where $\{z\}=\Ru{P(\mu)}\cdot{x}\cap Z(\mu)$, and $(K_{\mu})_x=(K_{\mu})_z$ is a maximal compact subgroup in $L(\mu)_x=L(\mu)_z$.
\end{proposition}

\begin{proof}
We follow the ideas of Brion's proof of Theorem~\ref{t:mom}.

First consider the case where $\mu$ is a rational point. Then $\tilde\mu=n\mu\in\Wt_X+n\wt$ is a lattice point in $n\Pol_{\Xc}$ for sufficiently big~$n$. Replacing $V$ with $\Sym^nV$ and $\delta$ with $n\delta$ by means of the Veronese map, which dilates $\Mom$ by the factor~$n$, we may assume that $n=1$ and $\mu$ is the eigenweight of a $B$-semi-invariant section $s\in\Ho0\bigl(X,\OO(\delta)\bigr)=V^*$.

Consider the irreducible $G$-module $V(\mu)$ with highest weight $\mu$ and the highest weight vector $u_\mu\in V(\mu)$. Choose a $K$-invariant Hermitian inner product on $V(\mu)$ (which is unique up to a positive multiple) and define a $K$-invariant Hermitian inner product on $V\otimes V(\mu)$ by the formula $$(v_1\otimes u_1\mid v_2\otimes u_2)=(v_1|v_2)\cdot(u_1|u_2),\qquad v_i\in V,\ u_i\in V(\mu).$$
Consider the variety $$\XG=\Xc\times Y_\mu\simeq G\times^{P(\mu)}\Xc\simeq K\times^{K_\mu}\Xc,$$ where $$Y_\mu=Gy_\mu\simeq G/P(\mu)\simeq K/K_\mu$$ is the unique closed orbit in~$\PP\bigl(V(\mu)\bigr)$, $y_\mu=\langle u_\mu\rangle$, and embed it by the Segre map into $\PP\bigl(V\otimes V(\mu)\bigr)$. The corresponding moment map is
\begin{equation}\label{eq:mom.Segre}
\Mom_\mu(x,y)=\Mom(x)+\Mom_\mu(y)
\end{equation}
(where $\Mom$ and $\Mom_\mu$ are the moment maps for the respective projective spaces). It is easy to see that
\begin{equation}\label{eq:mom.flag}
\Mom_\mu(y_\mu)=\mu.
\end{equation}
It follows that $\Mom(x)=-\mu$ if and only if $\Mom_\mu(x,y_\mu)=0$. In this case, by Proposition~\ref{p:KN-crit}\eqref{i:closed}, the orbit $G(x,y_\mu)$ is closed in the set of semistable points $\PP\bigl(V\otimes V(\mu)\bigr)^\sst$. Since $\Xc$ is spherical, it follows that $G$ has a dense open orbit in $\XG$ and therefore the closed $G$-orbit in $\XG^\sst_\mu=\XG\cap\PP\bigl(V\otimes V(\mu)\bigr)^\sst$ is unique.

Semistability means that there exists a homogeneous invariant polynomial $F\in\CC[V\otimes V(\mu)]^G$ of positive degree such that $F(v\otimes u_\mu)\ne0$, where $\langle v\rangle=x$. Applying the Veronese map once again, we may assume that $F\in V^*\otimes V(\mu)^*$ is a $G$-invariant linear form. Then the decomposition of $F$ over a tensor basis consisting of $T$-weight vectors is of the shape $F=s\otimes u^*_{-\mu}+{}\cdots{}$, where $s\in V^*$ is a highest weight vector of weight~$\mu$, $u^*_{-\mu}\in V(\mu)^*$ is the lowest weight vector (of weight~$-\mu$), and $F(v\otimes u_\mu)=s(v)$.

Therefore, the point $(x,y_\mu)$ lies in the unique closed $G$-orbit in ${\XG_F=\XG^\sst_\mu}$ and $x$ belongs to the unique closed $P(\mu)$-orbit in $\Xc_s=\Xc(\mu)$, which intersects $Z(\mu)$ in the unique closed $L(\mu)$-orbit. By Proposition~\ref{p:KN-crit}\eqref{i:crit}, the points $(x',y')\in\XG$ such that $\Mom_\mu(x',y')=0$ constitute a single orbit $K(x,y_\mu)$ and the points $x'\in\Xc$ such that $\Mom(x')=-\mu$ constitute the orbit $K_\mu\cdot{x}$.

By Proposition~\ref{p:KN-crit}\eqref{i:stab}, the stabilizer $K_{(x,y_\mu)}=(K_\mu)_x$ is a maximal compact subgroup in the reductive group $$G_{(x,y_\mu)}=P(\mu)_x=L(\mu)_x$$ (The last equality holds, because the Hermitian conjugation sends $P(\mu)$ to the opposite parabolic subgroup $P(-\mu)$, which intersects $P(\mu)$ in $L(\mu)$, while the group $P(\mu)_x$ is self-adjoint, since it is the complexification of~$(K_\mu)_x$.) If $x=uz$, $u\in\Ru{P(\mu)}$, then $$u^{-1}\cdot L(\mu)_x\cdot u=u^{-1}\cdot P(\mu)_x\cdot u=P(\mu)_z=L(\mu)_z.$$ It easily follows from the Levi decomposition of $P(\mu)$ that $L(\mu)_x=L(\mu)_z$ commutes with~$u$. This proves the last statement of the Proposition.

If $\mu$ is not rational, then, as it is well known, the symplectic structure on the coadjoint orbit $K\cdot\mu\simeq K/K_\mu\simeq G/P(\mu)=Y_\mu$ given by the Kirillov--Kostant--Souriau form $$\KKS(\xi\mu,\eta\mu)=\langle\mu,[\xi,\eta]\rangle,\qquad\forall\xi,\eta\in\kk,$$ defines a non-algebraic $K$-invariant K\"ahler structure, and the respective moment map $\Mom_\mu$ still satisfies~\eqref{eq:mom.flag}, where $y_\mu\in Y_\mu$ is the base point with stabilizer~$P(\mu)$. Indeed, for lattice points $\mu\in\Ch{T}$ the Fubini--Study metric on $Y_\mu\subset\PP(V(\mu))$ satisfies the above claim, and in the general case, $\mu$~is a positive linear combination of lattice vectors lying in the interior of the same face of the Weyl chamber, so that the desired K\"ahler metric on $Y_\mu$ can be defined as the respective linear combination of the Fubini--Study metrics.

For the product K\"ahler structure on $\XG=\Xc\times Y_\mu$, the moment map $\Mom_\mu$ is given by the formula~\eqref{eq:mom.Segre}. Replacing Proposition~\ref{p:KN-crit} by the results from~\cite[7.5]{Kir}, we get in the same way as before that if $\Mom(x)=-\mu$, then $(x,y_\mu)$ belongs to the closed $G$-orbit in the semistable locus (in the K\"ahler sense), i.e., in the open stratum $\XG^\sst_\mu$ of the Morse stratification of $\XG$ defined by the norm-square of the moment map $\Mom_\mu$; see~\cite[\S6]{Kir}. We show that
\begin{equation}\label{eq:sst}
\XG^\sst_\mu=G\cdot\bigl(\Xc(\mu)\times\{y_\mu\}\bigr)\simeq G\times^{P(\mu)}\Xc(\mu)
\end{equation}
does not depend on a point $\mu$ in the relative interior of~$\Pol(\mu)$.

For rational~$\mu$, the equality \eqref{eq:sst} holds by the above. By \cite[6.2]{Kir}, $(x,y_\mu)\in\XG^\sst_\mu$ if and only if $$\Mom_\mu\bigl(\overline{G(x,y_\mu)}\bigr)=K\cdot\Mom_\mu\bigl(\overline{P(\mu)x}\times\{y_\mu\}\bigr)\ni0,$$ which is equivalent to $\Mom\bigl(\overline{P(\mu)x}\bigr)\ni-\mu$ (here the upper bar denotes the closure).

Suppose $\Phi(x)=-\mu$. It follows from \cite[Thm.\,2.1]{conv} that $\Mom\bigl(\overline{P(\mu)x}\bigr)\cap\WC$ is a rational convex polytope contained in $-\Pol_{\Xc}$ and containing~$-\mu$. But then $\Mom\bigl(\overline{P(\mu)x}\bigr)$ also contains $-\mu'$ for some rational interior point $\mu'$ of~$\Pol(\mu)$. It follows that $(x,y_\mu)\in\XG^\sst_{\mu'}$, and $\XG^\sst_{\mu'}$ is given by~\eqref{eq:sst}. Thus $\XG^\sst_{\mu}\subseteq\XG^\sst_{\mu'}$.

It follows from the description of the unique closed $G$-orbit in $\XG^\sst_{\mu'}$ that if $\Mom(x')=-\mu'$, then $x'$ lies in the unique closed $P(\mu)$-orbit ${P(\mu)z\subseteq\Xc(\mu)}$. Since $\Mom$ is a proper map and $K_\mu$ acts transitively in the fibers of $\Mom$ over the relative interior of~$\Pol(\mu)$, for $\mu'$ sufficiently close to $\mu$ there exists $x'\in\Mom^{-1}(-\mu')$ arbitrarily close to~$x$. Hence ${x\in P(\mu)z}$, and $(x,y_\mu)$ lies in the unique closed orbit $G(z,y_\mu)\subseteq\XG^\sst_{\mu'}$. This yields the inverse inclusion ${\XG^\sst_{\mu}\supseteq\XG^\sst_{\mu'}}$ and finally the equality~\eqref{eq:sst}.

Besides, we have proved that $x$ belongs to the unique closed $P(\mu)$-orbit in $\Xc(\mu)$ whenever $\Mom(x)=-\mu$. By the Slice Theorem, for rational $\mu'\in\Pol(\mu)$ close to $\mu$, the stabilizer $(K_\mu)_{x'}$ of $x'\in P(\mu)x$ with $\Mom(x')=-\mu'$ is contained in a subgroup conjugate to~$(K_\mu)_x$. But by the above, $(K_\mu)_{x'}$ is a maximal compact subgroup in~$P(\mu)_{x'}=L(\mu)_{x'}$, and $(K_\mu)_x$ is equally so.
\end{proof}

\subsection{}

Now we describe the stratification of the orbit space $X/K$ in terms of the homeomorphism~\eqref{eq:X/K}.

\begin{theorem}\label{t:strat}
For any face $\Pol(\mu)$ of $\Pol_{\Xc}$ which does not coincide with any of the faces~$\Pol_{\Yc}$, where $Y\subset\Xc$ is a boundary $G$-orbit, the relative interior of $\Pol(\mu)$ is contained in one of the strata of~$X/K$. The respective orbit type is given by the conjugacy class of a maximal compact subgroup in the stabilizer $L(\mu)_z$ of the unique closed orbit $L(\mu)z$ in the BLV-slice $Z(\mu)$ of the affine chart
\begin{equation}\label{eq:B-chart}
X(\mu)=X\setminus\bigcup_iD_i\simeq P(\mu)\times^{L(\mu)}Z(\mu)\simeq\Ru{P(\mu)}\times Z(\mu),
\end{equation}
where the union is taken over all the divisors $D_i\subset X\setminus\Xo$ such that the inequalities \eqref{eq:pol} are strict.
\end{theorem}

\begin{proof}
Since $$\Kir^{-1}(-\mu)=\Mom^{-1}\bigl(K\cdot(-\mu)\bigr)=Kx\simeq K\times^{K_\mu}(K_\mu\cdot{x}),$$ where $K_\mu\cdot{x}=\Mom^{-1}(-\mu)$, the orbit type of $Kx$ is determined by the conjugacy class of the stabilizer~$(K_\mu)_x$. By Proposition~\ref{p:fib.mom}, $K_\mu\cdot{x}$~belongs to the unique closed $P(\mu)$-orbit in the affine chart~$\Xc(\mu)$, and $(K_\mu)_x$ is a maximal compact subgroup in~$L(\mu)z$, where $z$ is the unique intersection point of $\Ru{P(\mu)}$ and~$Z(\mu)$, so that $L(\mu)z$ is the unique closed orbit in~$Z(\mu)$.

It remains to observe that for any $G$-stable divisor $D_i\subset\Xc$ the inequality \eqref{eq:pol} turns into an equality exactly on the faces $\Pol_{\Yc}$ for the $G$-orbits $Y\subseteq D_i$, while these inequalities are strict on the relative interior of~$\Pol(\mu)$. Hence $\Xc(\mu)=X(\mu)$. (In~\eqref{eq:B-chart}, we alter the previous notation: $D_i$~denote rather not $B$-stable divisors on $\Xc$ intersecting~$X$, but their intersections with~$X$.)
\end{proof}

\begin{remark}
By appending back the removed faces $\Pol_{\Yc}$, we describe the stratification of the orbit space $\Xc/K$ in a similar way.
\end{remark}

\begin{remark}
As examples show (see~\S\ref{s:ex}), in general the stratification of the orbit space may be coarser than the stratification of the moment polytope by faces: a stratum in $X/K$ may correspond to a union of interiors of several faces of~$\Pol_{\Xc}$.
\end{remark}

\begin{remark}
Although Theorem~\ref{t:coarse} identifies the orbit space $X/K$ with the valuation cone~$\Val_X$, this correspondence remains mysterious in a sense. In particular, it is not clear how the faces of $\Val_X$ correspond to the faces of~$\Pol_{\Xc}$. We conjecture that there exists a natural homeomorphism between $\Val_X$ and $X/K$ which allows one to identify $\Val_X$ with ${\Pol_{\Xc}\setminus\bigcup_Y\Pol_{\Yc}}$ in an explicit way.

The strong version of Batyrev's conjecture \cite[Conj.\,6]{amoebae} says that the type of the $K$-orbit corresponding to a point in the relative interior of a face $\Face$ of $\Val_X$ is given by the conjugacy class of a maximal compact subgroup in the satellite subgroup $H_{\Face}\subseteq G$. All examples which we know show that the stabilizers for the $K$-orbits in $X$ and the maximal compact subgroups in the satellite subgroups constitute one and the same class of subgroups, up to conjugation. But, since maximal compact subgroups in $H_{\Face}$ may coincide for different faces $\Face\subseteq\Val_X$, as examples show, this does not allow us to associate a face of $\Val_X$ to each face of ${\Pol_{\Xc}\setminus\bigcup_Y\Pol_{\Yc}}$ anyway.
\end{remark}

Theorem~\ref{t:strat} gives some evidences in favor of the strong conjecture of Batyrev.

\begin{corollary}\label{c:gen.sat}
If $\mu=0$ or, more generally, $\mu$ is orthogonal to all the roots of~$G$, then the stabilizer for the respective $K$-orbit $\Kir^{-1}(-\mu)$ is a maximal compact subgroup in the ``general'' satellite subgroup $H_0=H$.
\end{corollary}

\begin{proof}
In this case, $P(\mu)=L(\mu)=G$, $X(\mu)=Z(\mu)=X$, and ${L(\mu)_z=H}$.
\end{proof}

\begin{corollary}\label{c:hor.sat}
The stabilizer for the $K$-orbit corresponding to an interior point of~$\Pol_{\Xc}$ is a maximal compact subgroup in the ``most degenerate'' satellite subgroup~$H_{\Val_X}$.
\end{corollary}

\begin{proof}
If $\mu$ is an interior point of~$\Pol_{\Xc}$, then all the inequalities \eqref{eq:pol} are strict and $X(\mu)=\Xo$ is the open $B$-orbit in~$X$. In this case, $L(\mu)$~acts on $Z(\mu)$ transitively and the derived subgroup $[L(\mu),L(\mu)]$ acts trivially. Furthermore, the closure of the BLV-slice $Z(\mu)$ in the open subset $$\Xco=\Xc\setminus\bigcup_{D_i\cap X\ne\varnothing}D_i\supseteq\Xo$$ is a smooth toric variety $Z$ with the fan~$\Fan_{\Xc}$, which intersects each $G$-orbit $Y\subseteq\Xc$ transversally in a single $L$-orbit, and $$\Xco=P(\mu)Z\simeq P(\mu)\times^{L(\mu)}Z\simeq\Ru{P(\mu)}\times Z;$$ see, e.g., \cite[Prop.\,3.3.2]{sph} or \cite[Thm.\,29.1]{hom-emb}. In particular, if $Y$ is closed, then $Y\cap Z=\{y\}$ is a single point with stabilizer $G_y=P(-\mu)$.

One easily deduces from \cite{sat} that $H_{\Val_X}={L(\mu)_z\ltimes\Ru{P(-\mu)}}$, where $z\in Z(\mu)$. Indeed, consider the blow-up $\Xbl\to X$ of the closed $G$-orbit $Y\subset X$. The vector $\bar\val_D$ corresponding to the exceptional divisor $D\subset\Xbl$ defines the barycentric subdivision of the cone $\Cone_Y$ and lies in the interior of~$\Val_X$. The normal bundle $N_{D,\Xbl}$ is equivariantly birationally isomorphic to $$N_{Y,\Xc}\simeq G\times^{P(-\mu)}T_yZ.$$ The space $T_yZ$ of dimension $n=\rk\Wt_X$ contains $P(-\mu)$-stable hyperplanes in general position $T_y(Z\cap D_{i_1}),\dots,T_y(Z\cap D_{i_n})$, where $D_{i_1},\dots,D_{i_n}$ are the $G$-stable divisors containing~$Y$. Hence $T_yZ$ decomposes into a direct sum of $P(-\mu)$-stable one-dimensional subspaces, on which $\Ru{P(-\mu)}$ acts trivially. The open $P(-\mu)$-orbit in $T_yZ$ is $L(\mu)$-equivariantly isomorphic to~$Z(\mu)$ and the stabilizer coincides with $L(\mu)_z\ltimes\Ru{P(-\mu)}$. This proves the claim.

Now, since a maximal compact subgroup of $H_{\Val_X}$ contained in $K$ is self-ad\-joint, it is a maximal compact subgroup in~$L(\mu)_z$. By Theorem~\ref{t:strat}, it coincides with the stabilizer for a $K$-orbit $\Kir^{-1}(-\mu)$ in the open stratum.
\end{proof}

\section{Examples}
\label{s:ex}

\begin{example}\label{e:flag}
Suppose that $X$ is projective, i.e., $X=G/P$, where $P\subset G$ is a parabolic subgroup. Such homogeneous varieties, including Grassmannians and classical varieties of complete or partial flags, are called \emph{generalized flag varieties}. Here is yet another characterisation: generalized flag varieties are homogeneous (under the full automorphism group) rational projective varieties.

It follows from the Bruhat decomposition in $G$ that a generalized flag variety $X$ contains a dense open orbit of the maximal unipotent subgroup~$\Ru{B}$, namely the big Schubert cell. Any $B$-semi-invariant rational function is constant on this orbit, whence $\Wt_X=0$, and also $\Val_X=\ES_X=0$. On the other side, it is well known that $K$ acts on $X$ transitively and the stabilizer of the base point $K_{x_0}=K\cap P$ is a maximal compact subgroup in~$P$.
\end{example}

\begin{example}\label{e:tor}
Now consider the case where $G=B=T=(\CC^\times)^n$ is an algebraic torus. In this case, $H$~coincides with the kernel of the action of $G$ on~$X$. Passing to the quotient by this kernel, we may assume that $H=\{\1\}$. The maximal compact subgroup $K=U_1^n$ is a compact $n$-di\-men\-sional torus (here $U_1\subset\CC^\times$ is the unit circle). The tangent Lie algebras are $\ttt=\CC^n$ and $\kk=\ii\RR^n$.

In this case, $\Wt_X=\Ch{T}\simeq\ZZ^n$ and $\Val_X=\ES_X\simeq\ii\kk=\RR^n$. The polar decomposition $$T=K\times\exp(\ii\kk)=U_1^n\times\RR_{>0}^n$$ immediately yields a homeomorphism $X/K\simeq\Val_X$. Here we have a unique orbit type, with trivial stabilizer.

On the other side, toroidal compactifications $\Xc$ considered in \ref{ss:mom} are in this case nothing but smooth projective toric varieties. If $\Xc\subseteq\PP(V)$ for some rational $T$-module~$V$ and a vector $v\in V$ lying over $x\in\Xc$ decomposes into a sum $v=v_1+\dots+v_m$ of $T$-eigenvectors $v_i$ of weights~$\lambda_i\in\Ch{T}$, then $$\Mom(x)=\Kir(x)=\sum_{i=1}^{m}\frac{(v_i|v_i)}{(v|v)}\lambda_i.$$ In particular, taking the base point $x=x_0\in X$ and an overlying vector $v=v_0$, we get
\begin{equation}\label{eq:mom.tor}
\Mom(tx_0)=\frac{\sum_ie^{2\lambda_i(\xi)}(v_i|v_i)\cdot\lambda_i}{\sum_ie^{2\lambda_i(\xi)}(v_i|v_i)},
\quad\forall t\in T,\ t=k\cdot\exp\xi,\ k\in K,\ \xi\in\ii\kk.
\end{equation}
The moment polytope $\Pol_{\Xc}\subset\ES_X^*$ is the convex hull of the eigenweights $-\lambda_1,\dots,-\lambda_m$ of $T$-semi-invariant linear functions not vanishing at~$v_0$. The formula \eqref{eq:mom.tor} in view of \eqref{eq:X/K} yields a one-to-one parametrization of the interior of $\Pol_{\Xc}$ by points in the dual vector space~$\ES_X$, which resembles the Fourier or Laplace transform. This parametrization should be known to experts in convex geometry, but the author did not succeed to find references in the literature.
\end{example}

\begin{example}\label{e:hor}
The case, where $X$ is a horospherical homogeneous variety, i.e., $H$~contains a maximal unipotent subgroup of~$G$, both generalizes and combines Examples \ref{e:flag} and \ref{e:tor} in a sense.

In this case, $P=N_G(H)$ is a parabolic subgroup in~$G$, with the Levi decomposition $P=L\ltimes\Ru{P}$, and $H=L_0\ltimes\Ru{P}$, where $L\supseteq L_0\supseteq[L,L]$ \cite[Lem.\,7.4]{hom-emb}. The variety $X$ has a structure of an equivariant principal bundle over the flag variety~$G/P$:
\begin{equation}\label{eq:hor}
X\simeq G\times^PA,
\end{equation}
where $A=P/H=L/L_0$ is an algebraic torus. The open $B$-orbit $\Xo$ is the preimage of the open $B$-orbit in~$G/P$.

Without loss of generality, we may assume that $L\supseteq T$ and $P$ contains the Borel subgroup opposite to~$B$. Then $x_0\in\Xo$, $A=T/(T\cap H)$, $\Wt_X=\Wt_A$, and $\Val_X=\ES_X^*=\ES_A^*$, where $A$ is regarded as a $T$-homogeneous variety \cite[7.2, 28.1]{hom-emb}.

It follows from the formula \eqref{eq:hor} and Examples \ref{e:flag}, \ref{e:tor} that $$X/K\simeq A/(T\cap K)\simeq\ES_A^*=\Val_X.$$ All $K$-orbits in $X$ have one and the same type, and stabilizers are conjugate to the maximal compact subgroup $K\cap L_0$ of~$H$.

Toroidal $G$-varieties with open orbit $X$ are of the shape $$\Xc=G\times^PZ,$$ where $Z$ is a toric $A$-variety \cite[5.3]{sph}. The variety $\Xc$ is smooth or projective if and only if $Z$ is so. The fan of the toric variety $Z$ equals~$\Fan_{\Xc}$ and the moment polytope for the embedding $Z\subset\Xc\hookrightarrow\PP(V)$ coincides with~$\Pol_{\Xc}$. A homeomorphism between $X/K$ and the interior of $\Pol_{\Xc}$ stems from Example~\ref{e:tor}.
\end{example}

\begin{example}\label{e:grp}
Consider the group variety $X=H$, where $H$ is a connected reductive algebraic group, on which the doubled group $G=H\times H$ acts transitively by multiplication from the left and from the right. In this case, $X=G/H$, where $H$ is embedded in $G$ diagonally, as the stabilizer of the base point $x_0=\1$. Sphericity of $X$ stems from the Bruhat decomposition in~$H$. Let us describe its invariants; see, e.g., \cite[27.2]{hom-emb}, \cite[5.4]{sph}.

For a Borel subgroup in~$G$, it is convenient to take $B=B_H^-\times B_H$, where $B_H^{\pm}$ are two mutually opposite Borel subgroups in~$H$, which intersect in a maximal torus~$T_H$. Choose a maximal torus $T=T_H\times T_H$ in~$B$. For a maximal compact subgroup in~$G$, we take $K=K_H\times K_H$, where $K_H$ is a maximal compact subgroup in $H$ compatible with~$T_H$.

The weight lattice $\Wt_X$ is identified with the character lattice $\Ch{T_H}$ embedded into $\Ch{T}$ antidiagonally: $\lambda\mapsto(-\lambda,\lambda)$. The valuation cone $\Val_X$ is identified with the negative Weyl chamber $$\WC_H\subset\ES_H=\Hom(\Ch{T_H},\RR)\simeq\ii(\ttt_H\cap\kk_H),$$ so that $G$-invariant valuations are determined via~\eqref{eq:val} by multiplicative one-parameter subgroups $x(t)\in H\bigl(\CC[t^{\pm1}]\bigr)$ corresponding to vectors in ${\WC_H\cap\Ch{T_H}^*}$.

The Cartan decomposition $$H=K_H\cdot\exp\ii\kk_H=K_H\cdot\exp\ii(\ttt_H\cap\kk_H)\cdot K_H$$ allows one to represent any $h\in H$ in the form $h=k_1\cdot\exp\val\cdot k_2$, where $k_1,k_2\in K_H$, and the point $\val\in\WC_H$ representing a $K$-orbit in $X$ is uniquely defined. This yields a homeomorphism $X/K\simeq\WC_H$.

For a vector $\val\in\WC_H$, the corresponding satellite subgroup is $$H_{\val}=\diag{L_H(\val)}\ltimes\bigl(\Ru{P_H(-\val)}\times\Ru{P_H(\val)}\bigr),$$ where $P_H(\val)$ and $L_H(\val)$ are a standard parabolic subgroup in $H$ and its Levi subgroup, respectively, given by~$\val$ \cite[1.1]{sat}. The stabilizer in $K$ of the point $x=\exp\val\in X$ is
$$K_x=\diag\bigl(K_H\cap L_H(\val)\bigr).$$ It is a maximal compact subgroup in~$H_{\val}$.
\end{example}

\begin{example}\label{e:sym}
A generalization of the previous example is the class of \emph{symmetric homogeneous varieties} $X=G/H$. Here $H$ coincides, up to connected components, with the subgroup $G^\theta$ of fixed points of an involutive automorphism $\theta$ of~$G$. If $G=H\times H$ and $\theta$ interchanges the factors, then $G^\theta=\diag(H)$, and we arrive at the group variety $X=H$ of Example~\ref{e:grp}.

Symmetric varieties are spherical. They were examined by Vust \cite{emb(symm)} in the framework of the general theory of spherical varieties. We give an account of necessary results; see, e.g., \cite[3.4.3]{sph} or \cite[\S26]{hom-emb} for details.

Consider the decomposition $\g=\g^\theta\oplus\g^{-\theta}$ into the direct sum of eigenspaces for $\theta$ with eigenvalues~$\pm1$, so that $\g^\theta=\h$. One can choose a $\theta$-stable maximal torus $T$ such that $T^\theta$ has minimal possible dimension. Then $Bx_0=\Xo$, and the weight lattice $\Wt_X=\Ch{T/(T\cap H)}\subseteq\Ch{T}$ consists of the characters taking value $1$ on~${T\cap H}$. The lattice $\Wt_X$ embeds into the space $\bigl(\ttt^{-\theta}\bigr)^*$ and spans its real form $\ii\bigl(\ttt^{-\theta}\cap\kk\bigr)^*$.

The valuation cone $\Val_X\subset\ES_X\simeq\ii\bigl(\ttt^{-\theta}\cap\kk\bigr)$ is identified with the Weyl chamber $\WC\cap\ii\bigl(\ttt^{-\theta}\cap\kk\bigr)$ of the \emph{little root system} of the symmetric variety~$X$, which is a root system in $\ii\bigl(\ttt^{-\theta}\cap\kk\bigr)^*$ consisting of the restrictions of the non-$\theta$-invariant roots of $G$ to~$\ttt^{-\theta}$. The corresponding Weyl group is $W_X=N_{K\cap H}(\ttt^{-\theta})/Z_{K\cap H}(\ttt^{-\theta})$ (the \emph{little Weyl group} of a symmetric variety).

There is a generalized Cartan decomposition \cite[Thm.\,4.1]{sym}:
$$G=K\cdot\exp\ii\bigl(\ttt^{-\theta}\cap\kk\bigr)\cdot H.$$
Moreover, for any $g\in G$, the decomposition $g=k\cdot\exp\val\cdot h$ with ${k\in K}$, ${h\in H}$, $\val\in\ii\bigl(\ttt^{-\theta}\cap\kk\bigr)$ can by adjusted so that $\val$ can be put into ${\WC\cap\ii\bigl(\ttt^{-\theta}\cap\kk\bigr)}$ by conjugation with an element of $N_{K\cap H}(\ttt^{-\theta})$, and such a vector $\val$ is uniquely determined by~$g$. This yields a homeomorphism $X/K\simeq\Val_X$.

A description of satellite subgroups in this example can be extracted from \cite[\S5]{comp.symm}. One can show that here the stabilizers $K_x$ of points $x=\exp\val$ are maximal compact subgroups in the respective satellite subgroups $H_\val$ as well.
\end{example}

Note that in Examples \ref{e:grp} and~\ref{e:sym}, the orbit type stratification of $X/K$ coincides with the face stratification of~$\Val_X$. This property breaks in the next example.

\begin{example}\label{e:SL_2^3/SL_2}
Let $G=\SL_2(\CC)\times\SL_2(\CC)\times\SL_2(\CC)$ and $H=\diag\SL_2(\CC)$. It is easy to see that the homogeneous variety $X=G/H$ is spherical: indeed, $B$-orbits in $X$ are in a one-two-one correspondence with $H$-orbits in $G/B=\PP^1\times\PP^1\times\PP^1$, the set of the latter orbits is finite and contains a dense open orbit. This example is considered in~\cite{amoebae}.

We identify $X$ with the closed subvariety in $\Mat_2(\CC)\times\Mat_2(\CC)\times\Mat_2(\CC)$ defined by the equations $x(1)x(2)x(3)=\1$, $\det{x(i)}=1$ ($i=1,2,3$) on a triple of matrices $x(1),x(2),x(3)\in\Mat_2(\CC)$. The group $G$ acts on $X$ by the rule: $$(g_1,g_2,g_3)\cdot\bigl(x(1),x(2),x(3)\bigr)=\bigl(g_2x(1)g_3^{-1},g_3x(2)g_1^{-1},g_1x(3)g_2^{-1}\bigr),$$
so that $(\1,\1,\1)\in X$ is the base point with stabilizer~$H$.

Choose a Borel subgroup $B$ consisting of the triples of upper-triangular matrices and a maximal torus $T$ consisting of the triples
$$t=\Biggl(
\begin{pmatrix}
t_1 &   0      \\
 0  & t_1^{-1} \\
\end{pmatrix},
\begin{pmatrix}
t_2 &   0      \\
 0  & t_2^{-1} \\
\end{pmatrix},
\begin{pmatrix}
t_3 &   0      \\
 0  & t_3^{-1} \\
\end{pmatrix}
\Biggr)$$
of diagonal matrices in~$\SL_2$. Take $K=\SU_2\times\SU_2\times\SU_2$ as a maximal compact subgroup in~$G$.

The lattice $\Ch{T}$ has a basis consisting of the fundamental weights $\omega_1,\omega_2,\omega_3$, where $\omega_i(t)=t_i$. The positive roots (which are also the simple roots) are $\alpha_i=2\omega_i$ (in the additive notation, i.e., $\alpha_i(t)=t_i^2$).

The $B$-stable prime divisors $D_{ij}\subset X$, $1\le i<j\le3$, bijectively correspond to the $H$-stable prime divisors on $\PP^1\times\PP^1\times\PP^1$ in the complement of the open $H$-orbit; the latter divisors consist of the triples $(p_1,p_2,p_3)$ of points in $\PP^1$ satisfying the condition $p_i=p_j$ for given $i,j$. In the coordinates on $X$, the divisor $D_{ij}$ is defined by the equation $f_{ij}(x):=x(k)_{21}=0$, where $\{i,j,k\}=\{1,2,3\}$.

The functions $f_{ij}$ generate the multiplicative group $\CC(X)^{(B)}/\CC^\times$ and their eigenweights $\omega_i+\omega_j$ generate the weight lattice~$\Wt_X$. The vectors in $\ES_X$ corresponding to the valuations $\val_{ij}=\val_{D_{ij}}$ are $\bar\val_{ij}=(\alpha_i^\vee+\alpha_j^\vee-\alpha_k^\vee)/2$, where $\alpha_1^\vee,\alpha_2^\vee,\alpha_3^\vee$ are the dual roots, which constitute a basis of $\Ch{T}^*$ dual to the basis of fundamental weights.

For each minimal parabolic subgroup $P_{\alpha_i}\subset G$, which consists of the triples $g=(g_1,g_2,g_3)$ of matrices in $\SL_2$ such that $g_j,g_k$ are upper-triangular, there are exactly two non-$P_{\alpha_i}$-stable divisors $D_{ij},D_{ik}$ in the complement of the open $B$-orbit. By a result of Luna (see, e.g., \cite[30.10]{hom-emb}), all the simple roots $\alpha_i$ are spherical roots. This means that the valuation cone $\Val_X$ coincides with the negative Weyl chamber~$\WC$.

For each pair $i,j$ of indices in $1,2,3$, we construct an equivariant compactification $\Xc_{ij}$ of~$X$. To this end, we identify $X$ with $\SL_2(\CC)\times\SL_2(\CC)\subset\Mat_2(\CC)\times\Mat_2(\CC)$ via the isomorphism $$\bigl(x(1),x(2),x(3)\bigr)\mapsto(x,y)=\bigl(x(i),x(j)^{-1}\bigr).$$  Suppose without loss of generality that $j$ follows $i$ in a cyclic permutation of $1,2,3$; then $G$ acts on $X$ by the following formula in the new coordinates:
$$(g_1,g_2,g_3)\cdot(x,y)=\bigl(g_jxg_k^{-1},g_iyg_k^{-1}\bigr).$$
Let us embed each of the two spaces $\Mat_2(\CC)$ into the projective space $\PP^4=\PP\bigl(\Mat_2(\CC)\oplus\CC\bigr)$ as an affine chart. Consider the matrix entries $x_{pq},y_{pq}$ as projective coordinates on~$\PP^4$, to which we also append the projective coordinates $x_0,y_0$ corresponding to the hyperplane at infinity. Define $\Xc_{ij}$ as the closure of $X$ in $\PP^4\times\PP^4$, given by the equations $\det{x}=x_0^2$, $\det{y}=y_0^2$. It is a smooth projective variety (a product of two smooth quadrics).

In projective coordinates on $\Xc_{ij}$, we have:
\begin{equation}\label{eq:fun}
f_{ij}=\frac{x_{22}y_{21}-x_{21}y_{22}}{x_0y_0},\qquad f_{ik}=-\frac{y_{21}}{y_0},\qquad f_{jk}=\frac{x_{21}}{x_0}.
\end{equation}
The numerators of these fractions are the equations of the closures of $D_{ij}$, $D_{ik}$, $D_{jk}$ in~$\Xc_{ij}$. The boundary $\Xc_{ij}\setminus X$ of the open orbit consists of the two $G$-stable prime divisors $D_i=\{y_0=0\}$ and $D_j=\{x_0=0\}$. It is easy to see from \eqref{eq:fun} that the valuations $\val_i=\val_{D_i}$ and $\val_j=\val_{D_j}$ correspond to the vectors $\bar\val_i=-\alpha_i^\vee$ and $\bar\val_j=-\alpha_j^\vee$ on two edges of~$\Val_X$.

The open $G$-orbits in these boundary divisors are $Y_i=D_i\setminus D_j$ and $Y_j=D_j\setminus D_i$. The intersection $D_i\cap D_j$ consists of the two orbits: the open one $Y_{ij}$ and the closed one~$Y_0$, which is defined by the condition $\Ker{x}=\Ker{y}$. Since $Y_0\subset\overline{D_{ij}}$, the variety $\Xc_{ij}$ is not toroidal. But there are no other $G$-orbits in $\Xc_{ij}$ contained in the closure of a $B$-stable divisor in~$X$. Therefore a toroidal smooth projective compactification of $X$ can be obtained by blowing up $Y_0$ and possibly some other $G$-orbits in the exceptional divisor.

The satellite subgroups corresponding to all faces of $\Val_X$, except for the vertex $\{0\}$ and the edge~$\RR_{\le0}\cdot\alpha_k^\vee$, can be computed as stabilizers for open $G$-orbits in the normal bundles of all orbits in~$\Xc_{ij}$, except~$Y_0$. Choose the following orbit representatives:
\begin{align*}
Y_i: &&
x&=
\begin{pmatrix}
1 & 0 \\
0 & 1 \\
\end{pmatrix},&
x_0&=1,&
y&=
\begin{pmatrix}
0 & 0 \\
0 & 1 \\
\end{pmatrix},&
y_0&=0;\\
Y_{ij}: &&
x&=
\begin{pmatrix}
1 & 0 \\
0 & 0 \\
\end{pmatrix},&
x_0&=0,&
y&=
\begin{pmatrix}
0 & 0 \\
0 & 1 \\
\end{pmatrix},&
y_0&=0.
\end{align*}
Their stabilizers in $G$ consist of the elements $g=(g_1,g_2,g_3)$ of the following shape:
\begin{align*}
Y_j: &&
g_i&=
\begin{pmatrix}
* & 0 \\
* & * \\
\end{pmatrix},&
g_j=g_k&=
\begin{pmatrix}
* & * \\
0 & * \\
\end{pmatrix};&&\\
Y_{ij}: &&
g_i&=
\begin{pmatrix}
* & 0 \\
* & * \\
\end{pmatrix},&
g_j&=
\begin{pmatrix}
* & * \\
0 & * \\
\end{pmatrix},&
g_k&=
\begin{pmatrix}
* & 0 \\
0 & * \\
\end{pmatrix}.
\end{align*}
These groups act in the normal spaces of the respective orbits at the above-indicated points. The satellites are the stabilizers of general points for these actions. They consist of the elements of the following shape:
\begin{align*}
Y_j: &&
g_i&=
\begin{pmatrix}
s &   0    \\
* & s^{-1} \\
\end{pmatrix},&
&g_j=g_k=
\begin{pmatrix}
s &  *     \\
0 & s^{-1} \\
\end{pmatrix};\\
Y_{ij}: &&
g_i&=
\begin{pmatrix}
s &   0    \\
* & s^{-1} \\
\end{pmatrix},&
&g_j=
\begin{pmatrix}
s &   *    \\
0 & s^{-1} \\
\end{pmatrix},\quad
g_k=
\begin{pmatrix}
s &  0     \\
0 & s^{-1} \\
\end{pmatrix}.
\end{align*}

On the other side, one may use Theorem~\ref{t:KM} to compute the orbit space~$X/K$. The normal space of the minimal $K$-orbit at the base point $x_0$ is $E\simeq\ii(\kk\cap\h)^\perp\subset\ii\kk$. It consists of the triples $(\xi_1,\xi_2,\xi_3)$ of matrices $\xi_i\in\ii\cdot\su_2$ satisfying the condition $\xi_1+\xi_2+\xi_3=0$. The group $K\cap H=\SU_2$ acts in $E$ diagonally by conjugation, and $X/K\simeq E/SU_2$.

Since the adjoint representation of $\SU_2$ in the space $\ii\cdot\su_2$ of traceless Hermitian matrices is equivalent to the standard linear representation of $\SO_3(\RR)$ in the Euclidean space~$\RR^3$, the space $E$ can be identified with the set of triples of vectors in $\RR^3$ with zero sum, on which $\SO_3(\RR)$ acts in a natural way. The invariants separating orbits are the vector lengths $$|\xi_i|=\sqrt{\tr\xi_i^2}=\sqrt2\cdot|x_i|,$$ where $\pm x_i$ are the eigenvalues of~$\xi_i$. They satisfy the triangle inequalities $$|\xi_i|\le|\xi_j|+|\xi_k|,\qquad\{i,j,k\}=\{1,2,3\},$$ which define a cone in $\RR^3$ homeomorphic to~$X/K$.

In the interior of this cone, all the inequalities are strict, the vectors $\xi_i$ span a plane in~$\RR^3$, the stabilizer of the triple $(\xi_1,\xi_2,\xi_3)$ in $\SO_3$ is trivial, and the stabilizer for the respective $K$-orbit in $X$ consists of $g=(g_1,g_2,g_3)$ such that
$$g_1=g_2=g_3=
\begin{pmatrix}
\pm1 &   0  \\
  0  & \pm1 \\
\end{pmatrix}.$$
On the boundary of the cone, except the vertex, where at least one of the inequalities is strict, the vectors $\xi_i$ are collinear, the stabilizer of the triple $(\xi_1,\xi_2,\xi_3)$ in $\SO_3$ is~$\SO_2$, and the stabilizer for the respective $K$-orbit in $X$ is conjugate to the subgroup consisting of the triples $g$ such that
$$g_1=g_2=g_3=
\begin{pmatrix}
s &   0    \\
0 & s^{-1} \\
\end{pmatrix},\qquad|s|=1.$$
At the vertex, $\xi_1=\xi_2=\xi_3=0$, the stabilizer is the whole group $\SO_3$, and the stabilizer for the respective $K$-orbit in $X$ is~$\diag\SU_2$.

In all cases, stabilizers of $K$-orbits are conjugate to maximal compact subgroups in satellite subgroups. But here the orbit type stratification of $X/K$ is coarser than the face stratification of~$\Val_X$.
\end{example}

\end{document}